\magnification=1200
\font\bb=msbm10
\def\square{\hbox{\rlap{$\sqcap$}$\sqcup$}}
\def\e{r_0^2/\kappa}
\def\d{r_0^4/\zeta}
\def\R{{\cal R}}
\def\z{\zeta}
\def\p{\psi}
\def\t{\theta}

\centerline{\bf $G$-manifolds with positive Ricci curvature and many isolated singular orbits}
\bigskip\bigskip
\centerline{\bf David J. Wraith}
\bigskip\bigskip
\itemitem{} {\bf Abstract:} {\it We show that in cohomogeneity 3 there are $G$-manifolds with any given number of isolated singular orbits and an invariant metric of positive Ricci curvature. We show that the corresponding result is also true in cohomogeneity 5 provided the number of singular orbits is even.}
\footnote{}{2000 {\it Mathematics Subject Classification:} 53C20, 53C21.}
\footnote{}{{\it Keywords:} $G$-manifold, cohomogeneity, positive Ricci curvature.}
\bigskip\bigskip
\centerline{\bf \S1 Introduction}
\bigskip\bigskip

The objects under consideration in this paper are compact $G$-manifolds with finitely many non-principal orbits. Here, $G$ is a compact Lie group acting smoothly and effectively on a smooth compact manifold $M.$ The orbits of such an action are either principal, exceptional (that is, non-principal but having the same dimension as a principal orbit), or singular, meaning that the orbit dimension is strictly lower than that of a principal orbit. The cohomogeneity of such an action is the codimension of a principal orbit, or equivalently the dimension of the space of orbits $G\backslash M.$ Note that the union of principal orbits is a dense subset of $M.$

We will be primarily interested in the invariant geometry of such objects, that is, the geometry of $M$ equipped with a Riemannian metric which is invariant under the $G$-action.

The motivation for studying manifolds with only finitely many non-principal orbits arose from the study of cohomogeneity one manifolds. Cohomogeneity one manifolds have been studied extensively in recent years, particularly for their geometric properties. (See for example the survey [Z].) Recently, a new example of a positively curved manifold was found among the cohomogeneity one manifolds ([D],[GVZ]), to add to the many known examples of non-negatively curved cohomogeneity-one manifolds (see for example [GZ1]). The underlying philosophy behind these developments is that curvature or other geometric expressions become simpler and more tractable in the prescence of symmetry, as provided by the action of a `large' Lie group acting isometrically.

Compact cohomogeneity one manifolds belong to one of two types according to the space of orbits. The orbit space could be a circle, in which case all orbits are principal and the manifold is a principal orbit bundle over the circle. The other possiblilty is a closed interval, in which case there are precisely two non-principal orbits corresponding to the end-points of the interval. This is by far the more interesting of the two cases, and provides our main motivating example for studying manifolds with a finite number of non-principal orbits. Such orbits are clearly isolated in the sense that there exist mutually disjoint invariant tubular neighbourhoods about each non-principal orbit. 

$G$-manifolds with finitely many non-principal orbits in cohomogeneities greater than one have been studied both topologically in [BW1] and geometrically in [BW2]. As it will be relevant later, we will first dicuss briefly the key topological features of these objects.

Let $K$ denote the principal isotropy of the $G$-action, and $H_1$,...,$H_p$ denote the non-principal isotropy groups. If $N_i$ denotes a tubular neighbourhood of the non-principal orbit $G/H_i,$ then $M^0:=M-\cup_{i=1}^p N_i$ has the structure of a principal-orbit bundle. Let $B$ denote the base of this bundle, so $B=G \backslash M^0.$ It is clear that $B$ is a manifold with $p$ boundary components. We note that $T_i:=\partial N_i$ has two fibration structures: it is fibered by principal orbits, and is also fibered by normal spheres $S^r$. The isotropy groups $H_i$ act on these normal spheres. If any $H_i$ acts transitively then the cohomogeneity must be one. As we are interested in cohomogeneities greater than one, we have that $H_i$ acts non-transitively, but with only one orbit type. It turns out that such actions are quite tightly constrained (see [B] chapter 4 \S6), and this results in the following
\proclaim Proposition 1.1. ([BW1], Theorem 9). If the cohomogeneity is greater than one, then $K$ is ineffective kernel of the $H_i$ action on $S^r,$ so $K$ is normal in $H_i$ and $H_i/K \cong U(1)$, $N_{SU(2)}U(1)$, $SU(2)$, or is finite, and acts freely and linearly on the normal sphere $S^r.$
\par

From this it is easy to deduce:
\proclaim Corollary 1.2. ([BW1], Corollary 10). If the cohomogeneity is greater than one, then $G\backslash T_i$ is either a complex or quaternionic projective space, or a $\hbox{\bb Z}_2$ quotient of an odd dimensional complex projective space in the case of a singular orbit, or in the case of an exceptional orbit a real projective or lens space. Also, each $G\backslash N_i$ is a cone over one of these spaces.
\par

In turn, we obtain the following description of the orbit space structure:
\proclaim Theorem 1.3. ([BW1], Theorem 3). $G\backslash M$ is the union of a manifold with boundary $B$, where each boundary component is one of the spaces listed in Corollary 1.2, together with a cone over each boundary component.
\par

Considering the dimensions of the possible boundary components of $B$ which can arise, we see that if our $G$-manifold contains a singular orbit, the cohomogeneity must be odd.
\medskip

It remains to describe the structure of a non-principal orbit neighbourhood. Let $L$ denote one of the groups listed in Proposition 1.1, and let $\alpha:L \rightarrow H_i/K$ be an isomorphism. There is a natural action of $L$ on the product $D^{r+1} \times G/K$ (where $D^{r+1}$ denotes a disc), given for $z \in L$ by $$z(x,gK) \mapsto (zx,g\alpha(z^{-1})K),$$ where the action of $L$ on $D^{r+1}$ is the standard Hopf action on the distance spheres about the origin of the disc. Using the symbol $\times_{\alpha}$ to denote a quotient under this action, we have
\proclaim Proposition 1.4. ([BW1], Theorem 3). For a small invariant tubular neighbourhood $N$ of a singular orbit $G/H,$ we have $N \cong D^{r+1} \times_{\alpha} G/K.$ For the tubular neighbourhood boundary $T$ we have $T \cong S^r \times_{\alpha} G/K.$
\par

We now turn our attention to the curvature of invariant Riemannian metrics. In the case of cohomogeneity one we have the following result:
\proclaim Theorem 1.5. ([GZ2]).  A compact cohomogeneity one manifold admits an invariant metric with positive Ricci curvature if and only if its fundamental group is finite.
\par

There is little possibility of proving a result as strong as this in the current context: the space of orbits in cohomogeneity one is either a circle or an interval. Either way, this makes no contribution to the curvature. However, in higher cohomogeneities, it is to be expected that the geometry of the space of orbits will play some role in determining the global geometric properties. 
\medskip

In [BW2] a general existence result is established for positive Ricci curvature metrics on $G$-manifolds with finitely many non-principal orbits (see Proposition 1.11 below). Using this result, many Ricci positive examples were presented. For instance we have the following collections, which feature the 7-dimensional Aloff-Wallach spaces as singular orbits. The Aloff-Wallach spaces are a $2$-parameter family of simply-connected $7$-dimensional homogeneous $\hbox{SU}(3)$-manifolds, which are very important in Riemannian geometry as almost all admit homogeneous metrics with positive sectional curvature (see [Z] page 82). Explicitly, for $p_1,p_2$ coprime, the Aloff-Wallach space $W_{p_1,p_2}$ is the quotient $$\hbox{SU}(3)/\Big\{\hbox{diag}(z^{p_1},z^{p_2},z^{-p_1-p_2}) \, | \, z\in\hbox{U}(1)\Big\}.$$ In [BW2] Theorems 5 and 6, it is shown that the families in Examples 1.6 and 1.7 below both contain infinitely many homotopy types.
\proclaim Examples 1.6. ([BW2]) Given any two Aloff-Wallach spaces $W_{p_1,p_2}$ and $W_{q_1,q_2}$, there is an 11-dimensional $\hbox{SU}(3)$-manifold $M^{11}_{p_1p_2q_1q_2}$ of cohomogeneity three, with orbit space $S^3$ and two singular orbits equal to the given Aloff-Wallach spaces, which admits an invariant metric of positive Ricci curvature.
\par

\proclaim Examples 1.7. ([BW2]) Given Aloff-Wallach spaces $W_{p_1,p_2}$ and $W_{q_1,q_2}$ with $p_1^2+p_1p_2+p_2^2=q_1^2+q_1q_2+q_2^2,$ there is a 13-dimensional $\hbox{SU}(3)$-manifold $M^{13}_{p_1p_2q_1q_2}$ of cohomogeneity 5, with orbit space a suspension of $\hbox{\bb C}P^2$ and two singular orbits equal to the given Aloff-Wallach manifolds, which admits an invariant metric of positive Ricci curvature. 
\par\medskip

The problem with the Ricci positive examples appearing in [BW2] is that all examples contain either one or two non-principal orbits. A natural question ([BW2] Open Problem number 1) is thus: is it possible to find invariant Ricci positive metrics on manifolds having more than two non-principal orbits? This is an intriguing question as the obvious candidates {\it just} fail. These candidates are those for which the manifold $B$ is a 3-sphere or a 5-sphere with some discs removed. Thus the boundary components are $\hbox{\bb C}P^1$s or $\hbox{\bb H}P^1$s respectively. It is easily checked that conditions (1) and (2) of Proposition 0.11 mean that while two discs can comfortably be removed, taking out three discs results in these conditions just failing to hold.
\par
The main aim of the current paper is to answer the above problem in the affirmative:
\proclaim Theorem 1.8. For any given $p \in \hbox{\bb N},$ there is a cohomogeneity three $\hbox{SU}(3)$-manifold with $p$ isolated singular orbits and an invariant metric of positive Ricci curvature, and a cohomogeneity five $\hbox{SU}(3)$-manifold with $2p$ singular orbits and an invariant metric of positive Ricci curvature.
\par
Theorem 1.8 follows immediately from the following theorem-examples:
\proclaim Theorem 1.9. Given any finite family of Aloff-Wallach spaces $\{W^7_1,...,W^7_p\},$ there exists an 11-dimensionsal $\hbox{SU}(3)$-manifold of cohomogeneity three with precisely $p$ singular orbits equal to the given Aloff-Wallach spaces, and an invariant metric of positive Ricci curvature.
\par
\proclaim Theorem 1.10. Given any finite family of Aloff-Wallach spaces $\{W^7_1,...,W^7_p\},$ there exists a 13-dimensionsal $\hbox{SU}(3)$-manifold of cohomogeneity five with precisely two singular orbits equal to each of the given Aloff-Wallach spaces, and an invariant metric of positive Ricci curvature.
\par
By an obvious adaptation of the proof of Theorem 6 in [BW1], we see that for each $p$, the families of manifolds in Theorems 1.9 and 1.10 contain infinitely many homotopy types.
\medskip

The construction behind the examples in Theorems 1.9 and 1.10 (which will be given in \S2) relies on two main ingredients, namely Propositions 1.11 and 1.12 below. In Proposition 1.11, $g_i$ denotes a metric on boundary component $i$ of $B$ which is induced via the standard submersion from the round metric of radius one. We adopt the convention that all principal curvatures at a boundary are computed with respect to the inward pointing normal.
\proclaim Proposition 1.11. ([BW2; Theorem 5]) Suppose that $\pi_1(G/K)$ is finite. Then $M$ admits an invariant Ricci positive metric if $B$ admits a Ricci positive metric such that \item{1)} the metric on boundary component $i$ is $\lambda^2_i g_i,$ and \item{2)} the principal curvatures (with respect to the inward normal) at boundary component $i$ are greater than $-1/\lambda_i.$ 
\par

To put this result in some context, notice that the principal part $M^0 \subset M$ clearly admits an invariant Ricci positive submersion metric if $B$ admits a Ricci positive metric, since the fibres of the submersion $G/K \hookrightarrow M^0 \rightarrow B$ admit Ricci positive metrics ([B; 9.70]). Moreover, a tubular neighbourhood of each non-principal orbit is easily seen to admit an invariant metric of non-negative sectional curvature. The issue addressed by Proposition 1.11 is essentially one of metric smoothness at the non-principal orbits: given a Ricci positive metric on $B,$ there is no guarantee in general that a corresponding submersion metric on $M^0$ can be smoothly extended over all of $M$ within Ricci positivity. Thus the point of Proposition 1.11 is to provide sufficient conditions under which this smoothing can be achieved.

\proclaim Proposition 1.12. Let $M$ denote the sphere $S^n,$ $n \ge 3,$ from which $p$ small non-intersecting discs have been removed. Then there is a Ricci positive metric on $M$ such that each boundary component is a round sphere of radius $\nu>0$ and all principal curvatures at the boundary are $>-1/(2\nu).$
\par

Combining Propositions 1.11 and 1.12 yields:
\proclaim Corollary 1.13. Suppose that $\pi_1(G/K)$ is finite, and that the space of orbits $B$ is either $S^3$ or $S^5$ with a number of non-intersecting discs removed. Then $M$ admits an invariant metric of positive Ricci curvature.
\par

\noindent{\bf Proof.} First note that $S^3$ or $S^5$ less a number of discs {\it is} a valid candidate for $B,$ as the boundary components are all equal to $S^2=\hbox{\bb C}P^1$ respectively $S^4=\hbox{\bb H}P^1.$ By Proposition 1.11, we only need to show that $B$ can be equipped with a Ricci positive metric satisfying conditions (1) and (2). By Proposition 1.12 we can equip $B$ in either case with a Ricci positive metric with round boundary components of radius $\nu$ with principal curvatures $>-1/(2\nu).$ Now the standard Fubini-Study metric on $\hbox{\bb C}P^1$ or $\hbox{\bb H}P^1$ is identical to a round metric of radius $1/2.$ Denoting the appropriate Fubini-Study metric by $g,$ we have that a round metric of radius $\nu$ is precisely $\lambda^2g$ with $\lambda=2\nu.$ Thus by Proposition 1.11, $M$ will admit an invariant Ricci positive metric provided the principal curvatures at the boundary are all $>-1/\lambda_i=-1/(2\nu),$ which is true by Proposition 1.12. \hfill\square

Using Corollary 1.13, we are now in a position to prove Theorems 1.9 and 1.10.
\medskip
\noindent {\bf Proof of Theorem 1.9.} Given Aloff-Wallach spaces $W_1,...,W_p,$ the $\hbox{SU}(3)$-manifold in question is just the obvious generalization to $p$ singular orbits of the manifold constructed in the proof of Theorem 5 in [BW1]. As the space of orbits for this manifold is $S^3$ with $p$ discs removed, the existence of an invariant Ricci positive metric follows immediately from Corollary 1.13.
\medskip

\noindent {\bf Proof of Theorem 1.10.} Given Aloff-Wallach spaces $W_1,...,W_p,$ consider the $\hbox{SU}(3)$-manifold $M_i^{13}$ which has two identical singular orbits equal to $W_i,$ as constructed in the proof of Theorem 26 in [BW1]. Away from the singular orbits, the $M_i$ have the structure of $\hbox{SU}(3)$-bundles. Performing fibre connected sums between the $M_i$ yields an $\hbox{SU}(3)$-manifold with singular orbits $W_1,W_1,...,W_p,W_p$ and orbit space an $S^5$ with $2p$ discs removed. The existance of an invariant Ricci positive metric now follows from Corollary 1.13. \hfill\square
\medskip

Theorem 1.8 leaves us with the following natural
\proclaim Open question. Can we find manifolds with more than two non-principal orbits and an invariant Ricci positive metric in cohomogeneities $\neq 3,5$? 
\par

A key feature of cohomogeneity 3 and 5 is that if we take any closed manifold of dimension 3 or 5 and remove a disc, the resulting boundary is a projective space, namely $\hbox{\bb C}P^1$ respectively $\hbox{\bb H}P^1.$ Thus a punctured 3-sphere or a punctured 5-sphere can be taken as the manifold $B,$ as in Corollary 1.13. In cohomogeneities $\neq 3,5$ removing a disc will not produce projective space boundaries. Thus we need to work harder to find candidates for $B.$ For example $\hbox{\bb H}P^{2k+1}$, $\hbox{\bb C}P^{2k+1}$ and $\hbox{\bb R}P^{2k+1}$ are boundaries, and so we can create manifolds with boundary (by a connected sum on the interior of the bounding manifolds) having any selection of these spaces as boundary components. To understand the topology, and especially the geometry of such objects presents a challenge. Although we believe the answer to the above question will be yes, we suspect that constructing examples to show this will prove difficult.
\medskip
The author would like to thank Stefan Bechtluft-Sachs and Fred Wilhelm for their valuable comments.
\bigskip\bigskip

\centerline{\bf \S2 Metrics on punctured spheres}
\bigskip\bigskip

In order to fully establish Theorems 1.8-1.10, it remains to prove Proposition 1.12 concerning Ricci positive metrics on punctured spheres. This is our objective in the current section.

In [P], Perelman provided a framework for constructing Ricci positive metrics on punctured spheres with round boundary components satisfying certain convexity conditions. In order to establish Proposition 1.12, we must show that within this framework we can exercise sufficient control to be able to achieve precisely those boundary conditions specified in the Proposition. In particular, this will involve establishing ranges and inter-relationships between various parameters needed in the construction.

Our task is made more difficult by the fact that in [P], constructions which are both subtle and complicated are for the most part only presented in outline. Consequently, in order to perform the desired analysis it is necessary to provide missing arguments for some of Perelman's claims - specifically those appearing in the sequel as Corollary 2.16, Lemma 2.17 and Lemma 2.18. In fact, we were not always able to provide such arguments using parameter values suggested by Perelman (the values of $\epsilon$ and $\delta$ given in [P; p162] proved problematic), and so alternative approaches have been developed as necessary. We believe that the ideas in [P] may be of wider interest, and for this reason some of the details provided in this section might be useful in other contexts.
\medskip 

Consider a round sphere $S^n,$ $n \ge 3,$ and remove a collection of non-intersecting discs. It is easy to check that the resulting boundaries do not satisfy the conclusion of Proposition 1.12. Moreover, it is also intuitively clear that we cannot smoothly glue tubes $S^{n-1} \times I$ to the boundary components so as to satisfy the boundary radius and principal curvature conditions whilst simultaneously maintaining positive Ricci curvature. Perelman's insight in [P] was that the kind of construction we want to achieve in Proposition 1.12 becomes easier if we `squash' the original sphere.

Consider a warped product metric $$dt^2+\cos^2 t\ ds^2_1 + R^2(t) ds^2_{n-2}$$ on $S^n,$ where $t \in [0,\pi/2].$ If $R(t)=\sin t$ then this metric is round of radius 1. The `squashing' suggested above is based on the singular metric $$dt^2+\cos^2 t\ ds^2_1 + R_0^2\sin^2 t\ ds^2_{n-2},$$ where $R_0<1$ is a constant. Given this metric, the sphere looks like a `flying saucer', with a singular circle corresponding to $t=0.$ The function $R(t)$ which will be constructed below will take the form $N\sin (t/N)$ for some $N>0$ in a small neighbourhood of $t=0,$ which ensures that the metric is smooth there. For larger values of $t,$ $R(t)$ will take the form $R_0\sin t.$ 

Next, we will remove a number of small geodesic discs of radius $r_0$ (see Definition 2.8 below) which are centred on the circle $t=0.$ The resulting boundary components are elliptical rather than round. To achieve the round boundaries needed for Proposition 1.12, we have to add tubes $S^{n-1} \times I$ to the boundary components, with the metric on the tubes chosen so as to interpolate from an elliptical boundary component to a round one. Of course, the tube metrics must also give the correct principal curvatures at the `outer' boundary, and glue smoothly with the punctured sphere at the `inner' boundary to give a globally Ricci positive metric.

It turns out that the tube metric construction and gluing results in \S2 respectively \S4 of [P] (see Propositions 2.20 and 2.21 of this paper) are sufficient for our purposes as stated without the need for further analysis. Thus our main task in establishing Proposition 1.12 is to perform a detailed construction of the punctured sphere metric.

For the ease of the reader we have tried to stay as close as possible to Perelman's notation, so as to facilitate comparisons between our computations and the constructions in [P].
\medskip

\proclaim Definition 2.1. Positive constants $R_0,$ $\kappa$ and $\z$ are chosen (in that order) according to the following rules:
\item{(a)} $R_0 \le 1/10;$
\item{(b)} $\kappa>2/\sqrt{3R_0};$
\item{(c)} $\z \in (\kappa, 3R_0\kappa^3/4).$
\par

Note that in the above definition, condition (b) ensures that the interval in condition (c) is non-empty. The choice of 1/10 as an upper bound for $R_0$ in condition (a) is somewhat arbitrary. We will need $R_0<1$ in the sequel, but difficulties arise if $R_0$ is too close to 1. Setting $R_0 \le 1/10$ means that such difficulties are easily avoided. 
\medskip

The following five numerical lemmas are needed for the definition of $r_0,$ the radius of the geodesic discs to be removed from $S^n.$

\proclaim Lemma 2.2. There exists a number $c_1=c_1(\kappa,\z)>0$ such that for all $0<x<c_1$ we have $${{\tan(x^2/\kappa+x^4/\zeta)} \over {\tan(x^2/\kappa)}}<1+\tan^2 x.$$ \par

\noindent{\bf Proof.} It is clear from the Taylor expansion of $\tan x$ that for $0<x<\pi/2$, $\tan x >x.$ Therefore, it suffices to show that $$\tan(x^2/\kappa+x^4/\zeta)<x^2/\kappa+x^4/\kappa,$$ for all $x$ sufficiently small. Expanding the left-hand side as a Taylor series, we obtain the inequality $$x^2/\kappa+x^4/\zeta + O(x^6) < x^2\kappa+x^4/\kappa,$$ which is true for sufficiently small $x$ provided $\z>\kappa.$ This in turn follows from Definition 2.1(c). \hfill\square
\medskip

\proclaim Lemma 2.3. There exists a number $c_2=c_2(\z)>0$ such that for all $0<x<c_2$ we have $$\sin(x+x^4/\z)\cot x <1.$$ \par

\noindent{\bf Proof.} We establish the equivalent inequality $\sin(x+x^4/\z)<\tan x$ for small $x.$ Expanding the left-hand side of this as a Taylor series we have $x-x^3/6+O(x^4),$ and expanding the right-hand side gives $x+x^3/3+O(x^5).$ Hence the inequality holds for all sufficiently small $x.$ \hfill\square
\medskip

\proclaim Lemma 2.4. There exists a number $c_3=c_3(\kappa)>0$ such that for all $0<x<c_3$ we have $$\tan^2 x > {{x^2} \over 2}+\tan^2({{x^2} \over {\kappa}}).$$
\par

\noindent{\bf Proof.} Using the fact that $\tan x>x$ for $0<x<\pi/2$, the inequality will certainly be true if $x^2/2>\tan^2(x^2/\kappa).$ Using the Taylor expansion for $\tan x$ we obtain $\tan^2(x^2/\kappa)=x^4/\kappa^2+O(x^8).$ Thus the desired inequality holds for all sufficiently small $x.$ \hfill\square
\medskip

\proclaim Lemma 2.5. There exists a number $c_4=c_4(R_0,\kappa,\z)>0$ such that for all $0<x<c_4$ we have $${1 \over {x^2}}\tan({{x^2} \over {\kappa}}+{{x^4} \over {\z}})>{{R_0} \over {\z}}\Bigl(1-{{x^3R_0} \over {\z}}\Bigr)^{-2}.$$
\par

\noindent{\bf Proof.} As $x \rightarrow 0,$ the left hand side tends to $1/\kappa,$ since $\tan x \approx x+x^3/3+O(x^5).$ The right hand side clearly tends to $R_0/\z.$ Thus the Lemma is proved if $1/\kappa>R_0/\z,$ that is if $\kappa<\z/R_0.$ But this is true by Definition 2.1(a) and (c).  \hfill\square
\medskip

\proclaim Lemma 2.6. There is a number $c_5=c_5(R_0,\kappa,\z)>0$ such that for all $0<x<c_5$ we have 
\item{(i)} $(x/2)\sin(2x/\kappa)>R_0\sin(x^2/\kappa+x^4/\z);$
\item{(ii)} $\cos(2x/\kappa)>R_0\cos(x^2/\kappa+x^4/\z);$
\item{(iii)} $$\cos(2x/\kappa)-R_0\cos(x^2/\kappa+x^4/\zeta)>{{(x/2)\sin(2x/\kappa)-R_0\sin(x^2/\kappa+x^4/\zeta)} \over {x^2/\kappa}}.$$
\par

\noindent{\bf Proof.} Expanding the various expressions in (i) and (ii) we have $$\eqalign{ {x \over 2}\sin(2x/\kappa)&={{x^2} \over {\kappa}}-{{2x^4} \over {3\kappa^3}}+O(x^6); \cr\cr
R_0\sin(x^2/\kappa+x^4/\z)&=R_0{{x^2} \over{\kappa}}+R_0{{x^4}\over {\z}}+O(x^6); \cr\cr
\cos(2x/\kappa)&=1-{{2x^2} \over {\kappa^2}}+{{2x^4} \over {3\kappa^4}}+O(x^6); \cr\cr
R_0\cos(x^2/\kappa+x^4/\zeta)&=R_0-R_0{{x^4} \over {2\kappa^2}}+O(x^6).\cr\cr}$$
Thus (i) and (ii) are clearly true if $x$ is small enough. For (iii), after multiplying both sides by $x^2/\kappa$ it suffices to establish the inequality $$(1-R_0){{x^2} \over {\kappa}}-{{2x^4} \over {\kappa^3}}+O(x^6)>(1-R_0){{x^2} \over {\kappa}}-x^4\Bigl({2 \over {3\kappa^3}}+{{R_0} \over {\z}}\Bigr)+O(x^6),$$ which for sufficiently small $x$ reduces to showing $${2 \over {\kappa^3}}<{2 \over {3\kappa^3}}+{{R_0} \over {\z}},$$ or equivalently $\z<3\kappa^3 R_0/4,$ and this is true by Definition 2.1(c). \hfill\square 

\proclaim Definition 2.7. Let $\gamma:\hbox{\bb R} \rightarrow \hbox{\bb R}$ be a standard smooth function with $\gamma' \le 0,$ interpolating between $\gamma(x)=1$ for $x \le 0$ and $\gamma(x)=0$ for $x \ge \Lambda,$ some $\Lambda>1/R_0,$ such that $$\sup_{x \in \hbox{\bb R}}\{|\gamma^{(k)}(x)|\} < R_0 \hbox{\ \  for } k=1,2.$$
\par
\noindent It is clear that such a function $\gamma$ exists for $\Lambda$ suitably large. 
\medskip

\proclaim Definition 2.8. Let $r_0>0$ be such that $r_0 < \min\{R_0,{{\pi} \over {2(1+\Lambda)}},c_1,c_2,c_3,c_4,c_5\}.$ \par
\medskip

In the next lemma we introduce a function $\R(t)$ which is a first approximation to our desired function $R(t).$ 

\proclaim Lemma 2.9. There exists $\p=\p(r_0)$ with $0<\p<r_0^2/\kappa$ and a $C^1$-function $\R(t)$ defined for $t \in [0,\pi/2]$ such that \item{(i)} $\R(t)=(r_0/2) \sin (2t/r_0)$ for $t \in [0,\p]$; \item{(ii)} $\R(t)=R_0 \sin(t+(\d)\gamma(t/r_0-1))$ for $t \in [\e,\pi/2]$; \item{(iii)} $\R(t)$ is smooth for $t \in (\p,\e)$ and $-\R''/\R \ge 4/r_0^2$ for these values of $t$. \par

For small $t,$ $\R(t)$ takes the form $N\sin (t/N)$ with $N=r_0/2,$ thus ensuring the smoothness of the metric at $t=0.$ For larger $t,$ we would like $\R(t)$ to be $R_0 \sin t.$ However, in order to achieve a $C^1$-join with the values specified for small $t$ (while maintaining the concavity requirement of 2.9(iii) which is necessary for curvature considerations) there has to be an adjustment, and this is achieved using the function $\gamma.$ This adjustment begins at $t =r_0>\e,$ and ends at $t=r_0(1+\Lambda).$ The condition $r_0<{{\pi} \over {2(1+\Lambda)}}$ in Definition 2.8 ensures that the adjusting effect of $\gamma(t)$ is exhausted by $t=\pi/2.$
\medskip

\noindent{\bf Proof of Lemma 2.9.} Since the form of $\R(t)$ is fixed for $t \in [0,\p]$ and $t \in [\e,\pi/2],$ we only have to show how to construct $\R(t)$ in the interval $t \in [\p,\e].$ Suitable values for $\p$ will emerge from the construction. 
\par
For $t \in [\p,\e]$ we consider the function $f(t):=(r_0/2)\sin(2t/r_0)+\t(t),$ for some function $\t(t) \le 0$ with $\t''(t) \le 0.$ To achieve a $C^1$-join between $f(t)$ and $(r_0/2)\sin(2t/r_0)$ at $t=\p$ we need $\t(\p)=\t'(\p)=0.$ To achieve a $C^1$-join at $t=\e$ we also need $$\eqalign{\t(\e)&=R_0\sin(\e+\d)-(r_0/2)\sin(2r_0/\kappa); \cr \t'(\e)&=R_0\cos(\e+\d)-\cos(2r_0/\kappa). \cr}$$ In order to complete the proof, we need to show that such a function $\t(t)$ exists, and that the resulting function $f(t)$ satisfies condition (iii) above.

First of all note that by Lemma 2.6 and the choice of $r_0,$ both $\t(\e)$ and $\t'(\e)$ must be negative. We can therefore clearly choose non-positive functions $\t(t)$ which satisfy the boundary conditions. The difficulty is doing this whilst preserving the concavity. By elementary calculus we can choose a concave-down $\t$ satisfying the boundary conditions if and only if $\t'(\e)$ is strictly more negative than the slope of the straight line joining the points $(\p,\t(\p))=(\p,0)$ and $(\e,\t(\e)),$ that is if $$\cos(2r_0/\kappa)-R_0\cos(\e+\d)>{{(r_0/2)\sin(2r_0/\kappa)-R_0\sin(\e+\d)} \over {\e-\p}}. \eqno{(\ast)}$$ Notice that if we establish the inequality $$\cos(2r_0/\kappa)-R_0\cos(\e+\d)>{{(r_0/2)\sin(2r_0/\kappa)-R_0\sin(\e+\d)} \over {\e}},$$ the openness of this inequality condition guarantees the existence of a very small $\p>0$ for which $(\ast)$ holds. However this second inequality is precisely the inequality appearing in Lemma 2.6(iii) evaluated at $x=r_0,$ and holds for our choice of $r_0$ by Lemma 2.6 and Definition 2.8. Thus a suitable value for $\p>0$ and a concave-down function $\t(t)$ satisfying the given boundary conditions can be found.

For $t \in [\p,\e]$ we now set $\R(t)=f(t),$ with $f(t)$ defined using the function $\t(t)$ chosen above. Thus $\R(t)$ is a $C^1$-function defined on the interval $[0,\pi/2].$ It remains to show that on $(\p,\e)$ we have $-\R''/\R\ge 4/r_0^2.$ However, note that since $\t(t)\le 0$ we have $\R(t) \le (r_0/2)\sin(2t/r_0)$ for these values of $t.$ Moreover, since $\t''(t) \le 0$ we also have $\R''(t) \le -(2/r_0)\sin(2t/r_0).$ Thus $$-{{\R''(t)} \over {\R(t)}} \ge {{(2/r_0)\sin(2t/r_0)} \over {(r_0/2)\sin(2t/r_0)}}.$$ As the right-hand side of the above expression is equal to $4/r_0^2,$ condition (iii) is established. \hfill\square
\medskip

\proclaim Lemma 2.10. For all $t \in (\e,\pi/2]$ we have $\R''(t)<0.$
\par

\noindent {\bf Proof.} For these values of $t$ we have $$\R(t)=R_0 \sin(t+(\d)\gamma(t/r_0-1)),$$ and so $$\eqalign{\R''(t)=&-R_0\sin(t+(r_0^4/\z)\gamma(t/r_0-1))(1+(r_0^3/\z)\gamma'(t/r_0-1))^2 \cr &+R_0\cos(t+(r_0^4/\z)\gamma(t/r_0-1))(r_0^2/\z)\gamma''(t/r_0-1). \cr}$$ As $\sup_{x \in \hbox{\bb R}}\{|\gamma^{(k)}(x)|\} < R_0\hbox{ for } k=1,2$ by definition of $\gamma,$ we see that $$\eqalign{\R''(t)<-&R_0\sin(t+(r_0^4/\z)\gamma(t/r_0-1))(1-(r^3_0R_0/\z))^2 \cr +& R_0\cos(t+(r_0^4/\z)\gamma(t/r_0-1))(r_0^2R_0/\z),\cr}$$ where we have replaced $\gamma''$ by $R_0$ in the second term, and the (non-positive) $\gamma'$ in the first term by $-R_0.$ Note that replacing $\gamma'$ in this way increases the value of the expression since $1-(r_0^3R_0/\z)>0.$ To see this last point, recall that by Definition 2.8, $r_0<R_0,$ and so this inequality will follow from the inequality $1>R_0^4/\z.$ But by Defintion 2.1 we have $\z>\kappa>2/\sqrt{3R_0},$ and so it suffices to show that $1>R_0^{9/2}{\sqrt 3}/2.$ Since $R_0 \le 1/10,$ this inequality is true.  

Therefore the result is established if we can show that $$\tan(t+(r_0^4/\z)\gamma(t/r_0-1))>(r_0^2R_0/\z)(1-(r_0^3R_0/\z))^{-2}.$$ As $\tan x$ is an increasing function, it suffices to consider the case where $t=r_0^2/\kappa,$ that is $${1 \over {r_0^2}}\tan(\e+\d)>{R_0 \over {\z}}\Bigl(1-{{r_0^3R_0} \over {\z}}\Bigr)^{-2}. \eqno{(\dag)}$$ Note that the smallest value of $t+(\d)\gamma(t/r_0-1)$ for $t \in [r_0^2/\kappa,\pi/2]$ {\it does} occur at $t=r_0^2/\kappa,$ since the reducing effect of $\gamma$ only begins at $t=r_0.$ For $t \ge r_0$ we have $t+(\d)\gamma(t/r_0-1)>r_0,$ so it suffices to show that $r_0>\e+\d,$ which follows from the fact that $r_0<R_0 \le 1/10$ and $\z>\kappa>2/\sqrt{3R_0}$ by Definition 2.1.

Finally, note that $(\dag)$ is true by Lemma 2.5 and the choice of $r_0.$
\hfill\square
\medskip

The following Lemma will be used in the proof of Lemma 2.15.
\proclaim Lemma 2.11. There exists $\iota=\iota(r_0)>0$ such that $$1-{{\R'} \over {\R}}(t) \tan t > \iota>0$$ for all $t \in [\e,r_0].$
\par

\noindent {\bf Proof.} For these values of $t,$ $\R(t)=R_0\sin(t+(r_0^4/\z)),$ and therefore $${{\R'} \over {\R}}(t) =\cot (t+(r_0^4/\z)),$$ which is strictly decreasing in $t.$ So $${{\R'} \over {\R}}(t) \tan t < \cot t \tan t =1.$$ By the compactness of the interval $[\e,r_0],$ the existence of $\iota$ follows. \hfill\square
\medskip

We now show how to smooth $\R(t).$
\proclaim Lemma 2.12. There exists a number $\mu_0=\mu_0(r_0,R_0,\kappa,\z)>0$ such that \item{(i)} $\mu_0<\p=\p(r_0,R_0,\kappa,\z);$ \item{(ii)} $\mu_0 \tan (r_0)<\iota=\iota(r_0),$ where $\iota$ is the quantity from Lemma 2.11; \item{(iii)} for all $t \in [\e,r_0],$ $$\mu_0<{{\cot t - \cot (t+(r_0^4/\z))} \over {1+\cot t}};$$ \item{(iv)} for all $t \in [\e,r_0],$ $$\mu_0<{{\cot(t+\d)(1+\cot^2 r_0)\tan t-\cot^2 r_0} \over {\tan r_0 (1+\cot^2 r_0)}}.$$
\par

\noindent {\bf Proof.} Conditions (i) and (ii) are easily fulfilled. Condition (iii) can be fulfilled as $\cot t$ is a strictly decreasing function of $t$ for $t \in (0,\pi/2).$ For (iv), it suffices to show that the numerator of the expression on the right-hand side is strictly positive, in other words $$\cot(t+\d)(1+\cot^2 r_0)\tan t>\cot^2 r_0.$$ This rearranges to $$\tan(t+\d)/\tan t < 1+\tan^2 r_0.$$ Computing the derivative of the left hand side shows that this quantity is strictly decreasing if and only if $\sin(2t)<\sin(2(t+\d)),$ which is true for the values of $t$ under consideration. Thus the maximum of the left hand side for $t \in [\e,r_0]$ occurs at $t=\e.$ Therefore the last inequality is true if it holds at $t=\e.$ But this follows from Lemma 2.2 and our choice of $r_0.$ \hfill\square
\medskip
Note that condition (ii) of Lemma 2.12 is used in Lemma 2.15, and (iii) and (iv) appear in Lemma 2.17.

\proclaim Lemma 2.13. Given any $\mu \in (0,\mu_0),$ we can smooth the function $\R(t)$ to a function $R(t)$ by adjusting the values of $\R(t)$ in the intervals $(\p-\mu,\p)$ and $(\e,\e+\mu),$ so that \item{(a)} $-R''/R > 2/r_0^2$ for all $t \in [0,\e]$; \item{(b)} $-R''/R >1-\mu$ for all $t \in [\e,\e+\mu]$; \item{(c)} $$\Big|{{R'(t)} \over {R(t)}}-{{\R'(t)} \over {\R(t)}}\Big|<\mu$$ for $t \in [\p-\mu,\p]$ and $t \in [\e,\e+\mu].$
\par

\noindent{\bf Proof.} We can smooth $\R(t)$ over the given intervals keeping both the values of the smoothed function and the values of its first derivative arbitrarily close to the original, and the second derivatives interpolating approximately linearly between the those on either side of the smoothing intervals. That conditions (a)-(c) can be satisfied by such a smoothing follows easily from the fact that $-\R''/\R \equiv 4/r_0^2$ when $t \in [0,\p),$ $-\R''/\R \ge 4/r_0^2$ when $t \in (\p, \e),$ and $-\R''/\R \equiv 1$ for $t\in(\e,r_0].$ \hfill\square
\medskip

\proclaim Corollary 2.14. For the smooth function $R(t)$, we have $R''(t)<0$ for all $t \in [0,\pi/2],$ and $-R''/R>1-\mu$ for all $t \in [\e,r_0].$
\par

\noindent{\bf Proof.} The first of these statements follows from Lemma 2.10, and (a) and (b) of Lemma 2.13. The second statement follows from 2.13(b), together with the observation that $-R''/R \equiv 1$ for $t \in (\e+\mu,r_0].$ \hfill\square
\medskip

\proclaim Lemma 2.15. For all $t \in [0,r_0]$ we have $${{R'} \over R} \tan t \le 1,$$ with the inequality being strict for $t \in (0,r_0].$
\par

\noindent {\bf Proof.} Using l'H\^opital's rule we see that $$\lim_{t \to 0^{+}} {{R'(t)} \over {R(t)}}\tan t=1.$$ For $t>0$ we need to check that $R'(t)\sin t < R(t)\cos t.$ As these are equal in the limit as $t \rightarrow 0^{+},$ it suffices to compare derivatives, and in particular the result will follow if we can establish $(R'(t)\sin t)' < (R(t)\cos t)'$ for $t>0.$
\par
Now $$\eqalign{(R'(t)\sin t)'=&\ \ \ R''(t)\sin t +R'(t)\cos t; \cr (R(t)\cos t)'=& -R(t)\sin t+R'(t)\cos t. \cr}$$ For $t \in [0,\e]$ we have $R'' << -R$ by 2.13(a), and thus the result follows for these values of $t.$
\par
For $t \in [\e,r_0],$ by Lemma 2.11 there exists $\iota>0$ such that $$1-{{\R'} \over {\R}}(t) \tan t > \iota.$$ By 2.13(c) we have $${{R'} \over R}(t)<{{\R'} \over {\R}}(t)+\mu.$$ Therefore the result will follow for these values of $t$ if $$\Bigl( {{\R'} \over {\R}}(t)+\mu \Bigr) \tan t <1,$$ that is, if $$\mu \tan t < 1-{{\R'} \over {\R}}(t) \tan t.$$ Thus it suffices to show that $\mu \tan (r_0)<\iota,$ and this inequality holds by our choice of $\mu_0.$ 
\item{} \hfill\square
\smallskip 

Next, we study principal curvatures at the boundary.
\proclaim Corollary 2.16. (Compare [P; \S3].) Equip $S^n,$ $n \ge 3$, with the metric $dt^2+\cos^2t\ ds_1^2+R^2(t)ds^2_{n-2}$ where $t \in [0,\pi/2].$ Remove a ball of radius $r_0$ centred on the circle $t=0.$ Then the principal curvatures at the resulting boundary are $\ge -\cot r_0.$
\par

\noindent{\bf Observation:} If $R(t)=\sin t$ for all $t \in [0,\pi/2]$ then the above metric is simply the unit radius round metric, and in this case the principal curvatures would all be identically equal to $-\cot r_0.$
\medskip

\noindent{\bf Proof.} A straightforward calculation of covariant derivatives shows that the principal curvatures occuring are $-\cot r_0$ and $- (R'(t)/R(t))\cot r_0 \tan t.$ Thus it suffices to show that $ (R'(t)/R(t))\tan t \le 1$ for $t \le r_0,$ and this is true by Lemma 2.15. \hfill\square
\medskip

In the next two lemmas, we investigate the sectional curvature of the intrinsic boundary metrics. Following Perelman, we denote the intrinsic sectional curvature by the symbol $K_i.$ For the rest of the notation, let $T=\partial/\partial t,$ and let $X$ denote a vector in the $S^1$-direction, with $T \wedge X$  denoting the plane spanned by these vectors. We will represent a vector tangent to $S^{n-2}$ by $\Sigma.$ Let $Y \in T \wedge X$ denote a vector tangent to the boundary. It might be helpful to note that the cosine of the angle between $T$ and the normal vector at any point on the boundary is $\cot r_0 \tan t.$ This follows from elementary spherical trigonometry.
\proclaim Lemma 2.17. (Compare [P; \S3].) The intrinsic curvatures $K_i(Y \wedge \Sigma)$ satisfy $$K_i(Y \wedge \Sigma) > \cot^2 r_0.$$
\par

\noindent{\bf Proof.} By [P] page 162, $K_i(Y \wedge \Sigma)$ is given by the expression $$K_i(Y \wedge \Sigma)=-{{R''} \over R}(1-\cot^2 r_0 \tan^2 t)+{{R'} \over R}\cot^2 r_0 \tan t(1+\tan^2 t).$$ This expression is not derived in [P], however it can be obtained by first computing the ambient sectional curvature using the formulas on page 159 of [P] (the latter formulas can themselves be obtained by computing Christoffel symbols, for example), then using the Gauss formula for the sectional curvature of embedded submanifolds (see [doC] page 130), and finally a little spherical trigonometry to obtain the form given above. 
\par
Consider $t \in [0,\e].$ For $t$ in this range we have $-R''/R > 2/r_0^2$ by Lemma 2.13(a), and thus $$K_i(Y \wedge \Sigma) > (2/r_0^2)(1-\cot^2 r_0 \tan^2 t).$$ As $\tan t$ is increasing with $t$, we see that $$K_i(Y \wedge \Sigma) > (2/r_0^2)(1-\cot^2 r_0 \tan^2(\e)).$$ Thus to show that $K_i(Y \wedge \Sigma)>\cot^2 r_0$ it suffices to show that $$(2/r_0^2)(1-\cot^2 r_0 \tan^2(\e))>\cot^2 r_0.$$ With a little rearrangement, this is equivalent to showing $$\tan^2 r_0 >(r_0^2/2)+\tan^2(\e).$$ But this is true by Lemma 2.4 and the choice of $r_0.$
\medskip
\noindent{\it Claim:} For all $t \in [\e,r_0],$ $K_i(Y \wedge \Sigma)>{{R'} \over R}(1+\cot^2 r_0)\tan t.$ 
\medskip
\noindent To establish this claim, recall from Corollary 2.14 that $-R''/R>1-\mu$ for $t$ in this range. Therefore $$K_i(Y \wedge \Sigma)>(1-\mu)(1-\cot^2 r_0 \tan^2 t)+{{R'} \over R}\cot^2 r_0 \tan t(1+\tan^2 t).$$ We therefore need to establish the inequality $$(1-\mu)-(1-\mu)\cot^2 r_0 \tan^2 t+{{R'} \over R}\cot^2 r_0 \tan t +{{R'} \over R}\cot^2 r_0 \tan^3 t \ge {{R'} \over R}(1+\cot^2 r_0)\tan t.$$ By gathering together the second and fourth terms on the left hand side, moving the third term on the left over to the right, and then simplifying the resulting inequality, we obtain $$\Bigl[(1-\mu)-{{R'} \over R}\tan t \Bigr] \ge \cot^2 r_0 \tan^2 t \Bigl[(1-\mu)-{{R'} \over R}\tan t \Bigr].$$ Assuming the term in the square brackets is non-negative, this inequality reduces to $1 \ge \cot^2 r_0 \tan^2 t,$ which is true since $t \le r_0$ by assumption.
\par
To complete the proof of the claim, it remains to show that $$(1-\mu)-{{R'} \over R}\tan t \ge 0.$$ For $t$ in the current range it follows from 2.13(c) that $$\cot(t+d)-\mu<{{R'} \over R}<\cot(t+\d)+\mu.$$ Therefore it suffices to show that $$1-\mu - \Bigl[\cot(t+\d)+\mu \Bigr] \tan t \ge 0.$$ This rearranges to $$\mu \le {{\cot t - \cot (t+(r_0^4/\z))} \over {1+\cot t}},$$ and this is true by our choice of $\mu_0.$ Thus the claim is established.
\par
Using the claim, to complete the proof of the Lemma it now suffices to show that $$(\cot(t+\d)-\mu)(1+\cot^2 r_0)\tan t> \cot^2 r_0$$ for $t \in [\e,r_0].$ This rearranges to $$\mu<{{\cot(t+\d)(1+\cot^2 r_0)\tan t-\cot^2 r_0} \over {\tan t (1+\cot^2 r_0)}}.$$ Since $\tan t \le \tan r_0$ for $t$ in the current range, it is enough to show that $$\mu<{{\cot(t+\d)(1+\cot^2 r_0)\tan t-\cot^2 r_0} \over {\tan r_0 (1+\cot^2 r_0)}}.$$ But this holds by the choice of $\mu_0$ in Lemma 2.12. \hfill\square
\medskip

In the next lemma, $\Sigma_1$ and $\Sigma_2$ are linearly independent tangent vectors to $S^{n-2}.$
\proclaim Lemma 2.18. (Compare [P; \S3].) The intrinsic curvatures $K_i(\Sigma_1 \wedge \Sigma_2)$ satisfy $$K_i(\Sigma_1 \wedge \Sigma_2) > \cot^2 r_0.$$
\par

\noindent{\bf Proof.} Perelman's claim ([P] page 162) that $$K_i(\Sigma_1 \wedge \Sigma_2)={{1-{R'}^2(t)(1-\cot^2 r_0 \tan^2 t)} \over {R^2(t)}}$$ is easily verified. We show (following Perelman) that $${{1-{R'}^2(t)(1-\cot^2 r_0 \tan^2 t)} \over {R^2(t)}} \ge {1 \over {\sin^2 t}}-\cot^2 t(1-\cot^2 r_0 \tan^2 t).$$ The right-hand side of this expression simplifies to $1+\cot^2 r_0,$ which is strictly greater than $\cot^2 r_0.$ Thus to establish the Lemma it suffices to establish the above inequality. 
\par
Notice that we would obtain equality in the above inequality if $R(t)=\sin t.$ Notice also that we can bound the left-hand side below by over-estimating both $R(t)$ and $R'(t)/R(t).$ We claim that for all $t \in [0,r_0],$ we have $$R(t) \le \sin t \hbox{ and } R'(t)/R(t) \le \cot t,$$ where the second of these statements follows immediately from Lemma 2.15. Thus establishing the first claim will complete the proof of the Lemma.
\par
For $t \in [0,\p-\mu]$ we need to check that $(r_0/2)\sin( 2t/r_0) \le \sin t.$ We have equality at $t=0,$ so comparing derivatives it suffices to show that $\cos (2t/r_0) \le \cos t,$ which requires $r_0 \le 2,$ and this is true by the choice of $r_0.$
\par
For $t \in [\p-\mu,r_0]$ we begin from the result (2.15) that $R'/R \le \cos t/\sin t,$ or equivalently $${d \over {dt}}\ln R \le {d \over {dt}} \ln (\sin t).$$
Integrating, we obtain $$\ln R(t)-\ln R(\p-\mu) \le \ln (\sin t)-\ln(\sin (\p-\mu)).$$ As $R(\p-\mu) \le \sin(\p-\mu)$ we have $$\ln R(t) \le \ln(\sin t)+\ln R(\p-\mu)-\ln(\sin (\p-\mu)) \le \ln(\sin t).$$ Thus $R(t) \le \sin t$ as required. \hfill\square
\medskip

\proclaim Proposition 2.19. Let $M$ denote the sphere $S^n$ from which $p$ small, non-intersecting dics have been removed. Then $M$ admits a Ricci positive metric such that all principal curvatures at each boundary component are $\ge -1,$ the induced metric on each boundary component can be expressed in the form $g=ds^2+B^2(s)ds^2_{n-2}$ where $s \in [0,\pi\omega],$ $1>\omega>\tau^{(n-2)/(n-1)}$ with $\tau:=\max B(s),$ and $g$ has all sectional curvatures $>1.$
\par

\noindent{\bf Proof.} Begin with the metric $dt^2+\cos^2 t\ ds_1^2 +R^2(t)ds^2_{n-2}$ on $S^n$ as above. Remove $p$ non-intersecting open balls of radius $r_0$ centered on the circle $t=0.$ Note that we are free to select a smaller value for $r_0$ (see Definition 2.8) should $p$ be too large for our original choice of $r_0.$ By Corollary 2.16, all principal curvatures at the boundary are $\ge -\cot r_0.$ By Lemmas 2.17 and 2.18 the sectional curvatures of the induced boundary metric are all $> \cot^2 r_0.$ Therefore, rescaling the metric on $M$ by a factor of $\cot^2 r_0$ produces principal curvatures $\ge -1$ and intrinsic sectional curvatures $>1,$ as required. It is clear that the (rescaled) metric on each boundary component can be expressed as $ds^2+B^2(s)ds^2_{n-2},$ where $B(s):=\cot(r_0)R(t(s)).$ We therefore have $$\tau:=\max B(s)=\cot(r_0)R(r_0)=R_0\cot(r_0)\sin(r_0+(\d)).$$ 

To find the range of the parameter $s,$ we must focus on the boundary component metric. In the first instance, we study the {\it unrescaled} metric.

We begin by observing that the metric $\sigma:=dt^2+\cos^2 t\ dx^2$ is a round unit radius metric on a hemisphere $D^2,$ where the boundary corresponds to $t=0.$ The metric on $S^n$ can then be viewed as a singular warped product metric $(D^2 \times S^{n-2}\,,\, \sigma+R^2(t)\ ds^2_{n-2}),$ with fibres collapsing at $t=0.$ Define $pr:S^n \rightarrow D^2$ to be the obvious projection map.

Let $x_0$ be a point on the boundary of $D^2,$ and remove an open (half) disc of radius $r_0$ centred on this point. Let $u(s)$ be a unit speed path along the arc $C$ created by removing this disc, so the length $L(u)$ is equal to the length of the arc. We claim that the preimage $pr^{-1}(C)$ is precisely the boundary of the distance sphere in $S^n$ of radius $r_0$ centered on the unique point in $S^n$ corresponding to $x_0 \in D^2.$ To see this, first note that topologically $pr^{-1}(C)$ is clearly a sphere of dimension $n-1.$ Furthermore, if $v(r)$ is any unit speed geodesic of length $r_0$ in $(D^2,\sigma)$ originating at $x_0,$ the `lifted' path $\hat{v}(r):=(v(r),y_0),$ where $y_0$ is any fixed choice of point in $S^{n-2},$ is a unit speed geodesic of length $r_0$ in $S^n.$ Thus $pr^{-1}(C)=\{(x,y)\,|\, x \in C \hbox{ and } y \in S^{n-2}\},$ which is the set of end-points of such geodesic lifts, must be a subset of the distance sphere $S(x_0,r_0) \subset S^n.$ But since the dimensions of these spheres are equal, the spheres themselves must coincide. 

Now consider the path $\hat{u}(s):=(u(s),y_0)$ in $S^n-D^n(x_0,r_0).$  By the above, this is a path in the boundary linking the two `poles' where the removed disc meets the circle $t=0.$ Moreover, $\hat{u}$ is a shortest path between these poles, as any path with non-constant $S^{n-2}$ coordinate must necessarily be strictly longer. Clearly $L(\hat{u})=L(u).$ 

We conclude that the distance between poles along the boundary created by removing a geodesic ball of radius $r_0$ from $S^n$ is equal to the length of the arc created by removing a ball of radius $r_0$ from $(D^2,\sigma)$ centered on the circle $t=0.$ By elementary spherical geometry, this length is easily seen to be $\pi \sin r_0.$

Now consider the rescaled metric. The corresponding distance between poles for this metric is $\pi\cot(r_0)\sin(r_0)=\pi\cos r_0,$ and hence the range of the parameter $s$ is $\pi\cos r_0.$ Therefore the constant $\omega$ in the statement of the Proposition takes the value $\cos r_0<1.$
\par
We need to check that $\omega>\tau^{(n-2)/(n-1)}.$ Firstly note that by Lemma 2.3 and the choice of $r_0,$ we have $\tau<R_0.$ As $\tau<1$ we have that $\tau^{(n-2)/(n-1)}$ is decreasing with $n.$ Thus it will suffice to show that $\omega>\sqrt{R_0}$ as $n \ge 3,$ or equivalently $r_0<\cos^{-1}(\sqrt{R_0}).$ As $r_0<R_0\le 1/10$ by definition of $r_0$ and $R_0,$ we see that this last inequality is true.
\par
It remains to show that the metric on $M$ has positive Ricci curvature. The metric on $S^n$ is essentially a warped product, as discussed above. The Ricci curvature formulas for such a metric are well-known (see for example [B] page 266): for the unrescaled metric $dt^2+\cos^2 t ds_1^2+R^2(t)ds^2_{n-2}$ we have $$\eqalign{\hbox{Ric}&(T,T)=1-(n-2)R''/R; \cr \hbox{Ric}&(X,X)=1+(n-2)(R'/R)\tan t; \cr \hbox{Ric}&(\Sigma,\Sigma)=-R''/R+ (R'/R)\tan t+(n-3)(1-{R'}^2)/R^2; \cr \hbox{Ric}&(T,X)=\hbox{Ric}(T,\Sigma)=\hbox{Ric}(X,\Sigma)=0. \cr}$$ Here, $T,$ $X$ and $\Sigma$ are as described before Lemma 2.17, but this time we assume in addition that all are {\it unit} vectors. Since $-R''(t)/R >0$ by Corollary 2.14 and $0 \le R' \le 1,$ we see immediately that this metric has positive Ricci curvature. Rescaling the metric simply rescales the Ricci curvature, and so has no effect on the positivity. \hfill\square
\medskip

\proclaim Proposition 2.20. ([P;page 159]) Let $g$ be a rotationally symmetric metric on $S^{n-1}$ with sectional curvature $>1,$ distance between the poles $\pi \omega$ and waist $2\pi \tau$; that is, $g$ can be expressed as $ds^2+B^2(s)ds^2_{n-2},$ where $s \in [0,\pi \omega]$ and $\max B(s)=\tau.$ Suppose that $\omega>\tau^{(n-2)/(n-1)},$ and let $\rho \in (\tau^{(n-2)/(n-1)},\omega).$ Then there exists a metric of positive Ricci curvature on $S^{n-1} \times [0,1]$ such that (a) the boundary component $S^{n-1} \times \{1\}$ has intrinsic metric $g$ and is strictly convex with all principal curvatures $>1$; (b) the boundary component $S^{n-1} \times \{0\}$ is concave with all principal curvatures equal to $-\lambda$ and is isometric to a round sphere of radius $\rho/\lambda$, for some $\lambda>0.$
\par

The idea is to glue a tube as described in Proposition 2.20 onto each of the boundary components of $M.$ To do this, we need the following gluing result:
\proclaim Proposition 2.21. ([P;\S4]) Suppose that $N_1$ and $N_2$ are compact smooth Riemannian manifolds with positive Ricci curvature and isometric boundaries. If the principal curvatures at $\partial N_1$ are strictly greater than the negatives of the corresponding principal curvatures at $\partial N_2,$ then the union $N_1 \cup N_2$ can be smoothed in a small neighbourhood of the gluing to produce a manifold of positive Ricci curvature.
\par

Combining Propositions 2.19, 2.20 and 2.21 in the obvious way, we arrive at the following:
\proclaim Corollary 2.22. For any choice of $\rho \in (\tau^{(n-2)/(n-1)},\omega),$ the manifold $M$ admits a Ricci positive metric such that each boundary component is a round sphere of radius $\rho/\lambda$ with all principal curvatures equal to $-\lambda.$ 
\par\medskip

\noindent{\bf Proof of Proposition 1.12.} We need to show that we can choose $\rho$ in Corollary 2.22 so that the boundary metrics are round of radius $\nu$ and the principal curvatures at the boundary are $>-1/(2\nu).$ From Corollary 2.22 we have $\nu=\rho/\lambda,$ and the principal curvatures all equal to $-\lambda.$ Thus the Proposition will be proved provided $-\lambda>-\lambda/2\rho,$ that is, provided $\rho<1/2.$ Now $\rho \in (\tau^{(n-2)/(n-1)},\omega).$ We have $\omega=\cos r_0,$ so this upper bound does not force $\rho<1/2.$ As $\rho$ can be taken to be any value in this interval, it therefore suffices to show that the lower bound $\tau^{(n-2)/(n-1)}<1/2.$ However, in the proof of Proposition 2.19 we argued that $\tau<R_0.$ Since $R_0 \le 1/10$ we have $$\tau^{(n-2)/(n-1)} \le \sqrt{\tau} < \sqrt{R_0} \le 1/\sqrt{10}.$$ As $1/\sqrt{10}<1/2$ the result follows. \hfill\square
\bigskip\bigskip

\centerline {\bf References}
\bigskip
\item{[B]} A.L. Besse, {\it Einstein Manifolds}, Springer-Verlag, Berlin (2002).
\medskip
\item{[BW1]}  S. Bechtluft-Sachs, D. J. Wraith, {\it On the topology of $G$-manifolds with finitely many non-principal orbits}, Topol. Appl. {\bf 159} (2012), 3282-3293.
\medskip
\item{[BW2]}  S. Bechtluft-Sachs, D. J. Wraith, {\it On the curvature of $G$-manifolds with finitely many non-principal orbits}, Geom. Dedicata {\bf 162} no. 1 (2013), 109-128.
\medskip
\item{[D]} O. Dearricott, {\it A 7-manifold with positive curvature}, Duke Math. J. {\bf 158} no. 2 (2011), 307-346.
\medskip
\item{[doC]} M. do Carmo, {\it Riemannian Geometry}, Birkh\"auser, (1992).
\medskip
\item{[GVZ]} K. Grove, L. Verdiani, W. Ziller, {\it An exotic $T_1S^4$ with positive curvature}, GAFA {\bf 21} no. 3 (2011), 499-524.
\medskip
\item{[GZ1]} K. Grove, W. Ziller, {\it Curvature and symmetry of Milnor spheres}, Ann. of Math. (2) {\bf 152} no. 1 (2000), 331-367.
\medskip
\item{[GZ2]} K. Grove, W. Ziller, {\it Cohomogeneity one manifolds with positive Ricci curvature}, Invent. Math. {\bf 149} (2002), 619-646.
\medskip
\item{[P]} G. Perelman, {\it Construction of manifolds of positive Ricci curvature with big volume growth and large Betti numbers}, `Comparison Geometry', Cambridge University Press, (1997).
\medskip
\item{[Z]} W. Ziller, {\it Examples of manifolds with non-negative sectional curvature}, Surveys in Differential Geometry volume XI, International Press (2007).

\bigskip\bigskip
\centerline{D. J. Wraith,}
\centerline{Department of Mathematics,}
\centerline{N.U.I. Maynooth,}
\centerline{Ireland.}
\smallskip
\centerline{email: David.Wraith@nuim.ie}

\end